\documentclass{amsart}
\usepackage{amssymb, amsmath, amsfonts, amscd}

\input xypic

\theoremstyle{plain}
\newtheorem{theorem}{Theorem}[section]
\newtheorem{corollary}[theorem]{Corollary}
\newtheorem{lemma}[theorem]{Lemma}
\newtheorem{proposition}[theorem]{Proposition}

\newtheorem{example}[theorem]{Example}

\theoremstyle{definition}
\newtheorem{definition}[theorem]{Definition}

\theoremstyle{remark}

\numberwithin{equation}{theorem}

\newcommand{\E}{\mathcal{E}}
\renewcommand{\O}{\mathcal{O} }

\renewcommand{\P}{\mathbf{P} }

\newcommand{\Hom}{\operatorname{Hom} }
\newcommand{\onabla}{\overline{\nabla} }
\newcommand{\End}{\operatorname{End} }
\newcommand{\Ext}{\operatorname{Ext} }
\newcommand{\diff}{\operatorname{diff}}
\newcommand{\Diff}{\operatorname{Diff} }
\renewcommand{\H}{\operatorname{H} }

\newcommand{\parone}{\partial_{x_1} }
\newcommand{\partwo}{\partial_{x_2} }
\newcommand{\parthree}{\partial_{x_3} }
\newcommand{\PARone}{\partial_{1} }
\newcommand{\PARtwo}{\partial_{2} }
\newcommand{\PARthree}{\partial_{3} }

\newcommand{\Spec}{\operatorname{Spec} }

\newcommand{\C}{\mathbf{C} }
\newcommand{\R}{\mathbf{R} }
\newcommand{\Der}{\operatorname{Der} }

\newcommand{\adn}{ad\nabla}

\newcommand{\pl}{\partial^l }

\newcommand{\Z}{\mathbf{Z}}

\begin{document}

\title{Higher order differential operators on projective modules}

\author{Helge Maakestad }
\address{H\oe gskolen i Bergen} 
\email{Helge.Maakestad@hib.no}



\date{October 2011} 

\begin{abstract}  The aim of this paper we introduce a new method - the "\emph{projective basis method}" -
  to give explicit formulas for generalized connections and higher order differential operators on any projective module.
We also give explicit formulas for logarithmic connections and connections on a  class of modules on ellipsiods.  Hence we get many explicit examples of modules
  on ellipsiods with non-flat connections and non-flat  higher order connections. The aim of the construction is to use the "projective basis -method" and higher order connections to  construct refined Chern classes in a refined algebraic DeRham cohomology. 
\end{abstract}

\maketitle
\tableofcontents

\section{Introduction}

For a finite rank locally trivial sheaf $E$  on a complex projective manifold $X$ there are few examples of non-flat connections
(see the the introduction of paper by Bloch and Esnault \cite{esnault}). 
If $E \cong \mathcal{O}_X^n$ is trivial of rank $n$, it follows $E$ has a flat connection defined using the universal derivation. The 
$\mathcal{O}_X$-module $E$ also has non-flat connections. The set of all connections on $E$ is parametrized by the finite dimensional vector space
\[ \Hom_{\mathcal{O}_X}(T_X, \End_{\mathcal{O}_X}(E)) \cong \Gamma(X, \Omega_X^1 \otimes \End_{\mathcal{O}_X}(E)),\]
and this vector space is non-trivial in general. Hence if $E$ has a connection, it has in general a large space of connections.  If $X$ is a complex projective manifold with non-trivial topological fundamental group $\pi_1(X)$
there is an equivalence of categories between the category pairs $(E, \nabla)$ consisting of a finite rank vector bundle $E$ on $X$ with a flat connection $\nabla$, and finite dimensional complex 
representations $(V, \rho)$ of the fundamental group $\pi_1(X)$. Hence to any non-trivial finite dimensional complex representation $(V, \rho)$ we get a non-trivial flat connection $(E(\rho), \nabla(\rho))$. Add a potential
$P\in \Gamma(X, \Omega_X^1 \otimes \End_{\mathcal{O}_X}(E))$ to get a new connection $\nabla:=\nabla(\rho) +P$. The new connection $\nabla$ is non-flat in general.  Hence the Riemann-Hilbert correspondence gives many examples of "non-trivial" non-flat connections. In Atiyah's paper  \cite{atiyah} it is conjectured that if a holomorphic vector  bundle on a complex projective manifold has a connection, then it has a flat connection. One of the subjects of this paper is to indicate that in the affine situation few vector bundles have flat connections. 

In the "differentiable category" for a smooth finite rank real vector bundle on a 
real smooth finite dimensional simply connected manifold, the existence of a flat connection implies that the bundle is a trivial bundle. 
This is mentioned without a proof on page 294 in Milnor's book \cite{milnor} in the appendix. This property does not hold in the algebraic category. In the paper \cite{maa1} we give examples of  non-free Cohen-Macaulay modules on isolated hypersurface singularities with flat algebraic connections. Hence flatness of the algebraic connection  does not imply triviality of the module. A trivial vector bundle on an affine algebraic variety always has a flat connection, hence if we  look for vector bundles on algebraic varieties with no flat algebraic connections we must consider nontrivial  vector bundles. 

One of the aims of this paper is to  give explicit formulas for non-flat algebraic connections on a large class of nontrivial finite rank vector bundles on ellipsoids  (see Theorem \ref{omega}) using a new method - the  "projective basis-method". This is a method giving explicit formulas for algebraic connections and higher order connections and differential operators on projective modules.
Another aim is to use projective bases to construct lifts of differential operators on the base variety to the vector 
bundle via higher order connections (see Definition \ref{canonicalconn}). 

In the "smooth category" it is known that any local operator 
is a differential operator (see \cite{peetre}), hence in analysis "most" operators are differential operators. In algebra any connetion is a differential operator of degree one. Higher order 
differential operators arise in the theory of $D$-modules and cristalline cohomology, hence differential operators are much studied in algebra, algebraic geometry and arithmetic geometry.

In a previous paper on a related subject the Kodaira-Spencer map and the Atiyah-class was used
to give explicit formulas for algebraic connections on maximal Cohen-Macaulay
modules on hypersurface singularities (see \cite{maa1}). In the papers
\cite{maa15} and \cite{maa2} the notion of an $I$-connection was
introduced where $I$ is a left and right $A$-module. The purpose of
this paper is to introduce the  "\emph{projective basis-method}" for
a finitely generated and projective $A$-module $E$ and to use such a
this method to give explicit formulas for $I$-connections on a class of
modules on ellipsiods (see Theorem \ref{basisconnection}). The connections we construct are non-flat in
general.

A projective basis for a finitely generated $A$-module $E$  (see the exercises in \cite{anderson}) is a set $e_1,\ldots, e_n$ of generators of $E$ and 
$x_1,\ldots, x_n$ of $E^*$ satisfying the relation
\[ \sum_i x_i(e)e_i=e \]
for all elements $e\in E$. An $A$-module $E$ has a projective basis if and only if it is finitely generated and projective.
Using a projective basis $B:=\{e_i,x_j\}$ for a finitely generated projective module $E$, we construct the fundamental matrix $\phi$ associated to 
$e_i, x_j$. This is an $n\times n$-matrix with coefficients in $A$ with $\phi^2=\phi$ - it is an idempotent endomorphism defininig the module  $E$. 
Using the projective basis $B$ we construct the "associated  connection $\nabla_B$"
\[ \nabla_B:\Der_k(A) \rightarrow \End_k(E) .\]

Two different projective bases $B,B'$ give rise to different connections $\nabla_B \neq \nabla_{B'}$. The set of connections on $E$ is a torsor on the $A$-module
$\Hom_A(\Der(A), \End_A(E))$: Given any "potential" $P \in \Hom_A(\Der(A), \End_A(E))$ it follows $\nabla:=\nabla_B + P$ is another connection on $E$. Hence the "space of connections" is parametrized by the $A$-module $\Hom_A(\Der(A), \End_A(E))$.

In Theorem \ref{curvature} we prove the following formula:

\begin{align}
\label{K1}  R_{\nabla_B}(\delta \wedge \eta)=[ \delta(\phi), \eta(\phi)]
\end{align}

for any derivations $\delta, \eta \in \Der_k(A)$. Hence the "curvature matrix" $R_{\nabla_B}(\delta \wedge \eta)$ is given as the Lie-product
of the two matrices $\delta(\phi), \eta(\phi)$ where $\phi$ is the idempotent matrix coming from the projective basis $B$ of $E$.  Here we view the curvature as an endomorphism

\[ R_{\nabla_B}(\delta \wedge \eta) \in \End_A(E).\]

The module $E$ is locally trivial, hence the curvature $R_{\nabla_B}(\delta\wedge \eta)$ is "locally a matrix". There is an open affine cover $U_i:=Spec(A_i)$ of $Spec(A)$ where the restriction $(R_{\nabla_B})_{U_i}$ is an $n\times n$-matrix with coefficients in $A_i$. In formula \ref{K1} we apply the derivations $\delta, \eta$ to the coefficients of the idempotent matrix $\phi$ and take the Lie product. This Lie product is non-trivial in general, hence the the connection $\nabla_B$ associated to a projective basis $B$ is non-flat in general. It follows the "projective basis-method" and formula \ref{K1} gives many examples of  non-flat algebraic connections. Since the space $\Hom_A(\Der(A),\End_A(E))$ is large in general, it is a non-trivial problem to determine if $E$ has a flat algebraic connection. One seeks a potential $P\in \Hom_A(\Der(A),\End_A(E))$ with the property that $\nabla:=\nabla_B+P$ is flat, and since $\nabla_B$ is non-flat in general this is a non-trivial problem.

In Example \ref{log} we give a construction of a "logarithmic connection"

\[ \nabla: \mathcal{E} \rightarrow \Omega^1_{X/k}(D) \otimes_{\mathcal{O}_X} \mathcal{E}, \]

where $X:=\Spec(A)$ and $\mathcal{E}:=\tilde{E}$ with $E$ a finitely generated  left $A$-module. Here $D:=V(f) \subseteq X$ is a hypersurface where $\mathcal{E}$ is locally trivial on the complement $U:=X-D$. The connection  $\nabla$ is constructed using a projective basis $\{B,B^*\}$ for the finite rank locally trivial $A_f$-module $E_f$. The connection
$\nabla$ is non-flat in general and the construction in this paper gives explicit formulas for the connection $\nabla$.

We also relate in Example \ref{BQconj} the construction to descent for modules as in \cite{nuss} and the Bass-Quillen conjecture. We give explicit examples of a finite rank projective 
module $E$ on $A[t]$ that is not extended from $A$ where $E$ has a flat connection relative to $A$. Hence for commutative rings there is no obvious relation 
between faithfully flat descent and the category of connections. Such a relation is described in \cite{nuss} for arbitrary associative rings, but Nuss considers a "non commutative connection" and a  version of the Amitsur complex.

We construct a higher order connection
\[ \nabla_B^l : E\rightarrow P^l \otimes_A E \]
associated to a projective basis $B$ for $E$ (see Proposition \ref{higher}). Here $P^l$ is the $l$'th
module of principal parts of the ring $A$.
There is a notion of flatness for higher order connections and
the connection $\nabla_B^l$ is non-flat in general.
We use the "projective basis-method" to construct a canonical map
\[ \rho_B:\Diff_\Z(A)\rightarrow \Diff_\Z(E) \]
of left $A$-modules (see Equation \ref{canonical}). 
The map $\rho_B$ is not a map of associative rings in general, hence the
left $A$-module structure on $E$ does not lift to a left
$\Diff_\Z(A)$-module structure on $E$. The obstruction to such a
lifting is given by the generalized curvature of the higher order connection $\nabla_B^l$.
To give a left  $\Diff_\Z(A)$-module structure on a left $A$-module $E$ is equivalent to giving a stratification on $E$ in the sense 
of cristalline cohomology (see \cite{ogus}).  The map $\rho_B$ is a ring-homomorphism if and only if $E$ is a left $\Diff_\Z(A)$-module. We show
that a canonical connection $\rho_B$ induced by a projective basis $B$ seldom is a ring-homomorphism, hence a finitely generated projective
module $E$ seldom has a stratification induced by a projective basis $B$. The notion of a stratification is used in the construction of cristalline 
cohomology (see \cite{ogus}).

\section{Differential operators and the "projective basis-method"}

Let in this section $A$ be a commutative ring with unit and let $E$ be
a left $A$-module. We give several criteria for the $A$-module $E$ to
be finitely generated and projective. We use the constructions to give
methods to calculate explicit $I$-connections and covariant
derivations on finitely generated projective modules. 

Note: The notion "projective basis" is in the litterature usually referred to as a "dual basis (see the exerciese in \cite{anderson}).
Some of the results on the relation between the existence of a projective basis and projectivity is considered to be  "well known", but I include it since I did not find a good reference.

Let in the following $B=\{e_1,..,e_n\}$ be a set of elements of $E$ and let
$B^*=\{x_1,..,x_n\}$ be a set of elements of $E^*$.

\begin{definition} We say the sets $\{B,B^*\}$ form a \emph{projective
    basis} for $E$ if the following holds:
\[ \sum_j x_j(e)e_j=e \]
for all elements $e\in E$.
\end{definition}


\begin{lemma} Assume $\{B,B^*\}$ is a projective basis for $E$. It
follows the set $B$ generates $E$ hence $E$ is finitely generated as
left $A$-module.
\end{lemma}
\begin{proof} Assume $e$ is an element in $E$. We get the equation
\[ e=\sum_j x_j(e)e_j \]
hence the set $B$ generates $E$ as left $A$-module since $x_j(e)\in A$
for all elements $e$ in $E$.
\end{proof}

Let $F=A\{u_1,..,u_n\}$ be a free left $A$-module on the basis
$C=\{u_1,..,u_n\}$, let $B=\{e_1,..,e_n\}$ be a set of elements of $E$
and define the following map:
\[ p: F\rightarrow E \]
by
\[ p(u_i)=e_i.\]
Let 
\[ y_i:F\rightarrow A \]
be defined by
\[ y_i(\sum a_ju_j)=a_i.\]
It follows $y_i=u_i^*$. 
Let $u=\sum_i a_iu_i\in F$. We get
\[ p(u)=\sum_i a_ip(u_i)=\sum_i a_i e_i=\sum_i y_i(u)e_i.\]

Define the following map
\[\rho :E^*\otimes_A E \rightarrow \End_A(E) \]
by
\[ \rho(\sum \phi_i \otimes e_i)(e)=\sum \phi_i(e)e_i.\]

\begin{definition} Let $end(E)=Coker(\rho)$ be $E's$ \emph{ring of projectivity}.
\end{definition}

Define the following product on $E^*\otimes_A E$:
\[ \bullet: E^*\otimes_A E\times E^*\otimes_A E\rightarrow E^*\otimes_A E \]
by
\[ \phi \otimes u \bullet \psi \otimes v=\psi\otimes \phi(v)u.\]

\begin{lemma} \label{endE}The following holds:
\begin{align}
&\label{ass1}\text{$\{ E^*\otimes_A E, \bullet\}$ is an associative ring.}\\
&\label{ass2}\text{$\rho$ is a map of $A$-algebras.}\\
&\label{ass3}\text{$Im(\rho)$ is a two sided ideal in $\End_A(E)$.}
\end{align}
\end{lemma}
\begin{proof}
Define the following map:
\[ f:E^*\times E\times E^* \times E \rightarrow E^*\otimes E \]
by
\[ f(\phi,u,\psi,v)=\psi\otimes\phi(v)u.\]
One checks $f$ is $A$-bilinear in all variables hence we get a well
defined product
\[\bullet: E^*\otimes E\times E^*\otimes E\rightarrow E^*\otimes E \]
defined by
\[ (\phi \otimes u)\bullet (\psi\otimes v)=\psi\otimes \phi(v)u.\]
One checks this product is left and right distributive over addition
in $E^*\otimes_A E$. Hence
\[ x\bullet(y+z)=x\bullet y + x\bullet z \]
and
\[ (y+z)\bullet x=y\bullet x + z\bullet x.\]
We check the product $\bullet$ is associative:
\[ \phi\otimes u\bullet(\psi \otimes v\bullet \chi \otimes w)\]
\[ \phi \otimes u\bullet \chi\otimes \psi(w)v=\]
\[ \chi\otimes \phi(\psi(w)v)u=\]
\[ \chi\otimes \psi(w)\phi(v)u=\]
\[ (\psi\otimes \phi(v)u)\bullet \chi \otimes w=\]
\[ (\phi \otimes u\bullet \psi \otimes v)\bullet \chi\otimes w.\]
It follows $\bullet$ is a product on $E^*\otimes_A E$ and hence
$E^*\otimes_A E$ is an
associative (nonunital) ring. One checks the map $\rho$ is a map of
rings and left $A$-modules hence Claim \ref{ass2} follows. We check
$Im(\rho)$ is a two sided ideal. It is clear $Im(\rho)$ is an abelian subgroup
of $\End_A(E)$. Assume $\psi\in \End_A(E)$ and $\phi\otimes e\in
Im(\rho)$.
It follows
\[ \psi\circ \phi\otimes e(x)=\psi(\phi(x)e)=\phi(x)\psi(e). \]
It follows
\[ \psi\circ \phi\otimes e=\rho(\phi\otimes \psi(e) ) \]
hence
\[ \psi\circ \phi\otimes e\in Im(\rho) \]
One checks
\[ \rho(\phi\otimes e)\circ \psi =\rho(\phi\circ \psi\otimes e) \]
hence $Im(\rho)$ is a two sided ideal in $\End_A(E)$
and Claim \ref{ass3} follows. The Lemma is proved.
\end{proof}

It follows from Lemma \ref{endE} $end(E)$ is an associative $A$-algebra.

Assume $\nabla: L\rightarrow \End_\Z(E)$ is a connection where
\[ \alpha:L\rightarrow \Der_\Z(A) \]
is a Lie-Rinehart algebra. This
means
$\nabla$ is an $A$-linear map and for all $a\in A$, $e\in E$ and
$\delta \in L$ it follows
\[ \nabla(\delta)(ae)=a\nabla(\delta)(e)+\delta(a)e.\]
The curvature of the connection $\nabla$ is the following map:
\[ R_\nabla: L\wedge_A L\rightarrow \End_A(E) \]
with
\[ R_\nabla(\delta\wedge \eta)=[\nabla(\delta),
\nabla(\eta)]-\nabla([\delta, \eta]).\]
The curvature $R_\nabla$ is the obstruction for $\nabla$ to be a map
of Lie algebras. We say $\nabla$ is \emph{flat} if $R_\nabla$ is zero.

Make the following definitions: The connection $\nabla$ induce a
canonical connection $\nabla^*$ on $E^*$ as follows: Let $\phi \in E^*$
and define
\[ \nabla^*(\delta)(\phi)=\delta\circ \phi-\phi\circ \nabla(\delta).\]
One verifies the map
\[ \nabla^*: L\rightarrow \End_\Z(E^*) \]
is a connection. We get a canonical connection
\[ \overline{\nabla}: L \rightarrow \End_\Z(E^*\otimes_A E) \]
defined by
\[ \onabla(\delta)(\phi\otimes e)=\nabla^*(\delta)(\phi)\otimes
e+\phi\otimes \nabla(\delta)(e).\]
There is a canonical connection
\[ ad\nabla: L \rightarrow \End_\Z(\End_A(E)) \]
defined by
\[ ad\nabla(\delta)(\phi)=\nabla(\delta)\circ \phi-\phi\circ
\nabla(\delta).\]

\begin{lemma} \label{map} The canonical map 
\[\rho:E^*\otimes_A E\rightarrow \End_A(E) \]
is a map of $L$-connections. The curvature of $ad\nabla$ is as
follows:
\[ R_{ad\nabla}(\delta\wedge \eta)(\phi)=[R_\nabla(\delta\wedge \eta), \phi].\]
\end{lemma}
\begin{proof}
We need to check the following: For any element $\phi \otimes e\in E^*\otimes_A E$ 
\[ ad\nabla(\delta)(\rho(\phi\otimes e))=\rho(\onabla(\delta)(\phi\otimes e)) .\]
We get:
\[ ad\nabla(\rho(\phi\otimes e))(x)=\nabla(\delta)(\phi(x)e)-\phi(\nabla(\delta)(x))e.\]
We see
\[ \rho(\onabla(\delta)(\phi\otimes e)(x)=\]
\[ \rho(\nabla^*(\delta)(\phi)\otimes e+\phi\otimes
\nabla(\delta)(e))(x)=\]
\[\nabla^*(\delta)(\phi)(x)e+\phi(x)\nabla(\delta)(e)=\]
\[\delta(\phi(x))e-\phi(\nabla(\delta)(x))e+\phi(x)\nabla(\delta)(e).\]
Since $\nabla$ is a connection it follows
\[
\delta(\phi(x))e+\phi(x)\nabla(\delta)(e)=\nabla(\delta)(\phi(x)e).\]
It follows
\[ ad\nabla(\rho(\phi\otimes e))=\rho(\onabla(\delta)(\phi\otimes
e)) \]
and the Lemma is proved. Assume $\delta, \eta \in L$. We get
\[ R_{ad\nabla}(\delta\wedge \eta)=[ad\nabla(\delta),
ad\nabla(\eta)]-ad\nabla([\delta, \eta]).\]
It follows
\[ R_{\adn}(\delta, \eta)(\phi)=\adn(\delta)\circ
\adn(\eta)(\phi)-\adn(\eta)\circ
\adn(\delta)(\phi)-\adn([\delta,\eta])(\phi)=\]
\[ \adn(\nabla(\eta)\circ \phi-\phi\circ
\nabla(\eta))-\adn(\eta)(\nabla(\delta)\circ \phi-\phi\circ
\nabla(\delta)) \]
\[ -(\nabla([\delta, \eta])\circ \phi -\phi \circ \nabla([\delta,
\eta])=\]
\[ [\nabla(\delta),\nabla(\eta)]\circ \phi -
\phi\circ[\nabla(\eta),\nabla(\delta)]-[\nabla([\delta,
\eta]), \phi]=\]
\[ [[\nabla(\delta),\nabla(\eta)],
\phi]-[\nabla([\delta,\eta]),\phi]=\]
\[ [ R_\nabla(\delta\wedge \eta), \phi].\]
It follows
\[ R_{\adn}(\delta\wedge \eta)(\phi)=[R_\nabla(\delta\wedge
\eta),\phi] \]
and the Lemma is proved.
\end{proof}

\begin{theorem} \label{projectiveequivalence} The following are equivalent:
\begin{align}
&\label{p1}\text{$E$ is a finitely generated projective $A$-module.} \\
&\label{p2}\text{$E$ has a projective basis.}\\
&\label{p3}\text{$id_E$ is in the image of $\rho$.}\\
&\label{p4}\text{$\rho$ is an isomorphism.}\\
&\label{p5}\text{$end(E)=0$}
\end{align}
Moreover if $E$ has an $L$-connection $\nabla$ it follows $end(E)$ has an $L$-connection
$\theta_\nabla$. If $\nabla$ is flat it follows $\theta_\nabla$ is flat.
\end{theorem}
\begin{proof} Define $x_j=y_j\circ s$ where $s$ is a section of the
  map $p:F\rightarrow E$. It follows $s(e_i)=u_i+y $ where $y$ is in
  $ker(p)$. Let $e=e_i$. We get
\[ \sum_j x_j(e)e_j=\sum_j y_j(s(e_i))e_j=\sum_j y_j(u_i+y)e_j=\]
\[\sum_j y_j(u_i)e_j+\sum_j y_j(y)e_j=\sum_j \delta_{ij}e_j+p(y)=e_i=e\]
since $y\in ker(p)$. It follows $\sum_j x_j(e_i)e_j=e_i$ for all
$i$. Assume
$e=\sum_i a_ie_i$. We get
\[  \sum_j x_j(e)e_j=\sum_j x_j(\sum_i a_ie_i)e_j=\]
\[ \sum_ia_i\sum_j x_j(e_i)e_j=\sum_ia_ie_i=e.\]
It follows the sets $\{B,B^*\}$ form a projective basis for
$E$. Conversely if $\{e_i,x_j\}$ is a projective basis for $E$ it
follows the map
\[ s(e)=\sum_j x_j(e)u_j\] 
is an $A$-linear section of $p$: 
\[ p(s(e))=p(\sum_j x_j(e)u_j)=\sum_j x_j(e)p(u_j)=\sum_j
x_j(e)e_j=e.\]
The equivalence of \ref{p1} and \ref{p2} is shown.
We prove the equivalence of \ref{p1} and \ref{p3}: Let $\omega=\sum_j
x_j\otimes e_j\in E^*\otimes E$ be an element. If $f(\omega)=id_E$ it
follows
\[ \sum_j x_j(e)e_j=id_E(e)=e \]
hence $\{B,B^*\}$ is a projective basis. It follows $E$ is finitely
generated and projective. Assume conversely $E$ is finitely generated
and projective with projective basis $\{B,B^*\}$. It follows 
\[ f(\sum_j x_j\otimes e_j)=id_E .\]
We have proved the equivalence of \ref{p1} and \ref{p3}. 
We prove the equivalence of \ref{p1} and \ref{p4}: Assume $E$ is
finitely generated and projective. It follows the map $f$ is an isomorphism. Assume
conversely $f$ is an isomorphism. It follows there is an element 
\[ \sum_j x_j \otimes e_j \in E^*\otimes E \]
mapping to $id_E$. The elements $x_j, e_i$ gives rise to a projective
basis $\{B,B^*\}$ for $E$ hence by the equivalence above it follows
$E$ is finitely generated and projective.
The equivalence between \ref{p1} and \ref{p5} is clear. 

Assume $\nabla$ is a connection on $E$. Let $ad\nabla$ be the induced
connection on $\End_A(E)$. By Lemma \ref{map} it follows the two sided ideal $Im(\rho)$ is
stable under the action of $ad\nabla$ hence $ad\nabla$ induces a
connection $\theta_\nabla$ on $end(E)$. By Lemma \ref{map} it follows
\[ R_{ad\nabla}(\delta\wedge \eta)(\phi)=[R_\nabla(\delta \wedge
  \eta),\phi].\]
Hence if $\nabla$ is flat it follows $ad\nabla$ is flat hence
the induced connection $\theta_\nabla$ is flat. The Theorem is proved.
\end{proof}

Hence the associative ring $end(E)$ is the obstruction for the module
$E$ to be finitely generated and projective.

Note: In the paper \cite{maa1} we used the Kodaira-Spencer map and
Atiyah class to calculate explicit expressions of flat connections 
on a class of maximal Cohen-Macaulay modules on surface
singularities. In this calculation we investigated the endomorphism
ring $\End_A(M)$ for a maximal Cohen-Macaulay module $M$ and a set of
special elements of this ring. It might be this calculation is related
to the ring of projectivity $end(M)$ of $M$. This topic will be
studied in future papers on the subject.

Consider the following example from \cite{maa1} Section 2: Let
$f=x^m+y^n+z^2$ be an element of $K[x,y,z]$ where $K$ is a field of
characteristic zero and let $A=K[x,y,z]/f$.
Let $\phi, \psi$ be the following matrices:
\[
\phi=\begin{pmatrix} x^{m-k} & y^{n-l} & 0 & z \\
                      y^l &  -x^k & z & 0 \\
                    z & 0 & -y^{n-l} & -x^k \\
                   0 & z  & x^{m-k} & -y^l 
\end{pmatrix}
\]
and
\[
\psi=
\begin{pmatrix} x^k & y^{n-l} & z & 0 \\
                y^l & -x^{m-k} & 0 & z \\
            0       & z & -y^l & x^k \\
            z &       0 & -x^{m-k} & -y^{n-l} 
\end{pmatrix}
\]
with $0\leq k \leq m-1$ and $0 \leq l \leq n-1$.
It follows the matrices $\phi$ and $\psi$ are a matrix factorization
of the polynomial $f$ and we get a periodic resolution
\[ \cdots \rightarrow ^\psi A^4 \rightarrow ^\phi A^4 \rightarrow^\psi A^4
\rightarrow^\phi M(\phi,\psi) \rightarrow 0 \]
where  $M=M(\phi,\psi)$ is a maximal Cohen-Macaulay module on $A$.
In Theorem 2.1 and 3.3 in \cite{maa1} we calculated a flat connection
\[ \nabla: \mathbb{V}_M \rightarrow \End_K(M) \]
for all pairs of matrices $\phi, \psi$ and all $m,n\geq 2$ and $k,l$. The module $M$ is not locally
free hence the ring of projectivity $end(M)$ is non-zero.

\begin{corollary} \label{flatm} There is for any pair of matrices $\phi,\psi$ a flat
  connection
\[ \theta: \mathbb{V}_M\rightarrow \End_K(end(M)) .\]
\end{corollary}
\begin{proof} The Corollary follows immediately from Theorem
  \ref{projectiveequivalence}, Theorem 2.1 and Theorem 3.3 in \cite{maa1}.
\end{proof}

We get from Corollary \ref{flatm} an $U(\mathbb{V}_M)$-module
structure on the associative ring $end(M)$ where $U(\mathbb{V}_M)$ is
the universal enveloping algebra of the Lie-Rinehart algebra $\mathbb{V}_M$.

Assume $E$ has a projective basis $\{B,B^*\}$.
Let $I$ be a left and right $A$-module and let $d\in
\Der_\mathbf{Z}(A,I)$ be a derivation.
By \cite{maa2}, Proposition 2.13 there is a characteristic class
\[ c_I(E)\in \Ext^1_A(E,I\otimes_A E) \]
with the property $c_I(E)=0$ if and only if $E$ has an $I$-connection
\[ \nabla: E\rightarrow I\otimes_A E \]
with
\[ \nabla(ae)=a\nabla(e)+d(a)\otimes e.\]
If $E$ is projective it follows $\Ext^1_A(E,I\otimes_A E)=0$ hence
trivially $c_I(E)=0$. We aim to give an explicit calculation of a
generalized connections in this case.
Define the following maps:
\[ \nabla: E\rightarrow I\otimes_A E \]
by
\[ \nabla(e)=\sum_i d(x_i(e))\otimes e_i .\]

Define also
\[ \nabla':\Der_\mathbf{Z}(A) \rightarrow \End_\mathbf{Z}(E) \]
by
\[ \nabla'(\delta)(e)=\sum_i \delta(x_i(e))e_i.\]

\begin{theorem} \label{basisconnection} The maps $\nabla, \nabla'$ are connections on $E$.
\end{theorem}
\begin{proof} Assume $a\in A$ and $e\in E$. We get
\[ \nabla(ae)=\sum_i  d(x_i(ae))\otimes e_i=\]
\[ \sum_i d(ax_i(e))\otimes e_i=\sum_i ad(x_i(e))\otimes e_i+\sum_i
d(a)x_i(e)\otimes e_i=\]
\[ a\sum_i d(x_i(e))\otimes e_i+ \sum_i d(a)\otimes x_i(e)e_i=\]
\[ a\nabla(e)+d(a)\otimes e.\]
We get moreover
\[ \nabla'(\delta)(ae)=\sum_i \delta(x_i(ae))e_i=\]
\[ \sum_i \delta(ax_i(e))e_i= \]
\[ \sum_i \delta(a)x_i(e)e_i+ a\delta(x_i(e))e_i=\]
\[a\nabla'(e)+\delta(a)e .\]
It follows the maps $\nabla, \nabla'$ are connections on $E$.
\end{proof}

We see from Theorem \ref{basisconnection} that the notion of a projective basis gives rise to explicit formulas for $I$-connections
where $I$ is any left and right $A$-module and $E$ is  any finitely generated projective $A$-module. 

\begin{example}\label{log}  Explicit formulas for logarithmic connections.\end{example}

In the study of "logaritmic differential forms" and "logarithmic connections"
one encounters the sheaf $\Omega^1_{X}(D)$ where $D \subseteq X:=\Spec(A)$ is a divisor with normal crossings. Let $(E,\nabla)$ be a finite rank projective $A$-module 
with a logarithmic connection $\nabla$. Theorem \ref{basisconnection} gives rise to explicit formulas for the connection $\nabla$. Hence the formalism has applications 
to the study of log geometry, the logarithmic deRham complex and logarithmic deRham cohomology. 

Let $k$ be a fixed commutative unital ring and let $A$ be a commutative $k$-algebra. If $E$ is a coherent $A$-module and if $f\in A$ and assume $D:=V(f) \subseteq X:=\Spec(A)$ is a hypersurface where $E_f$ is locally trivial on the complement $U:=D(f)$. Let $\mathcal{E}$ be the $\mathcal{O}_X$-module defined by $E$ and let $\mathcal{E}_f$ be the restriction of $\mathcal{E}$ to $U$. There is since $U$ is affine and $\mathcal{E}_f$ finite rank and locally trivial $\mathcal{O}_U$-module, a connection

\[  \nabla^*_B: \mathcal{E}_f \rightarrow \Omega^1_{U}\otimes_{\mathcal{O}_U} \mathcal{E}_f ,\]

where $\mathcal{E}_f$ is the sheafification of the localized module $E_f$ and $\{B,B^*\}$ is a projective basis for $E_f$. There is an isomorphism $\Omega^1_{A_f/k} \cong S^{-1}\Omega^1_{A/k}$ where $S:=\{f^n\}$ is the multiplicative set defined by powers $f^n$ of $f$. We get an isomorphism

\[ \Omega^1_{A_f/k} \otimes_{A_f} E_f \cong S^{-1}\Omega^1_{A/k}\otimes_A E \]

and an induced "connection" 

\[  \nabla: \mathcal{E} \rightarrow \Omega^1_{X/k}(D)\otimes_{\mathcal{O}_X} \mathcal{E} ,\]

where we have defined $\Omega^1_{X/k}(D)$ to be the sheafification of $S^{-1}\Omega^1_{A/k}$. Since $\nabla^*_B$ is non flat in general, it follows $\nabla$ is non-flat in general.

It follows Theorem \ref{basisconnection} gives explicit formulas for the "logarithmic connection" $\nabla$.  Hence if you choose an arbitrary finitely generated $A$-module $E$ and a projective basis $\{B,B^*\}$ for $E_f$, it follows the corresponding "logarithmic connection" $\nabla$ will in general be non-flat. If $D$ is a divisor with "normal crossings" and if $k$ is the field of complex numbers, the above construction is related to the notion "logarithmic connection" studied in complex analysis.

\begin{example} \label{BQconj}A relation with the the Bass-Quillen conjecture. \end{example}

Let $A$ be a commutative unital ring and let $B:=A[t]$ be the polynomial ring in the variable $t$ over $A$. Let $E$ be a finite rank projective $B$-module with 
projective basis $B$ and corresponding connection

\[ \nabla_B: \Der_{\mathbb{Z}}(B) \rightarrow \End_{\mathbb{Z}}(E).\]

There is an inclusion of Lie-Rinehart algebras $\Der_{A}(B) \subseteq \Der_{\mathbb{Z}}(B) $ and we get an induced "covariant derivative"

\[ \nabla: \Der_A(B) \rightarrow \End_{\mathbb{Z}}(E).\]

\begin{lemma} \label{BQ} The covariant derivative $\nabla$ is flat.
\end{lemma}
\begin{proof} Let $D:=\partial_t \in \Der_A(B)$ be partial derivative wrto the $t$ variable. It follows $\Der_A(B)$ is free on $D$ as a left $B$-module.
Let $x:=aD,y:=bD$ with $a,b \in B$ and consider the curvature 

\[  R_{\nabla}(x,y):= [\nabla(x), \nabla(y)]-\nabla([x,y]).\]

It follows

\[ [x,y]=(aD(b)-D(a)b)D \]

and 

\[  [\nabla(x),\nabla(y) ]= [a\nabla(D), b\nabla(D)]= (aD(b)-D(a)b)\nabla(D).\]

Hence

\[  R_{\nabla}(x,y)=  (aD(b)-D(a)b)\nabla(D)-\nabla((aD(b)-D(a)b) D)=0 \]
and the Lemma follows.
\end{proof}

Let $\Omega:=\Omega^1_{B/A}$ be the module of Kahler differentials of $B$ over $A$. It follows $\Omega \cong B dt$ is a free rank one $B$-module on the element $dt$,
hence $\wedge^2_B \Omega =0$ is trivially zero. Hence any connection

\[ \nabla: E \rightarrow \Omega \otimes_B E \]

is automatically flat.

In \cite{nuss} the notion \emph{flat connection} is related to descent for modules. If $R \rightarrow S$ is a faithfully flat map of commutative unital rings, one wants to classify left $S$-modules $M$ with the property that there is a
left $R$-module $N$ and an isomorphism of $S$-modules $S\otimes_R N \cong M$. This problem is related to the notion of a flat connection. The following results are proved:  Let $\operatorname{Mod}(R)$ be the category of left $R$-modules and let $\operatorname{Desc}(\psi)$ be the category of descent data for $\psi$.

\begin{theorem} \label{n1} There is an equivalence of categories  $\operatorname{Mod}(R) \cong  \operatorname{Desc}(\psi)$.
\end{theorem}
\begin{proof} See \cite{nuss} Theorem 3.8.
\end{proof}

Let $\operatorname{Flconn}(\psi)$ be the category of flat connections $(E, \nabla)$ as mentioned in section 3.4 in \cite{nuss} page 44, where 

\[ \nabla: E \rightarrow \Omega^1_{S/R}\otimes_S E \]

is a connection in the sense that $\nabla$ is an $R$-linear map with the property that for any $s\in S, e\in E$ it follows

\[  \nabla(se)=s\nabla(e)+ d(s)\otimes e.\]

\begin{theorem} Assume $\psi: R \rightarrow S$ is a faithfully flat map of commutative unital rings. \label{n2} There is an equivalence of categories $\operatorname{Desc}(\psi)  \cong \operatorname{Flconn}(\psi)$.
\end{theorem}
\begin{proof} See \cite{nuss}, Theorem 3.11.
\end{proof}

Note: By Traverso's paper \cite{traverso}, a reduced commutative unital ring $A$ has the following property: There is an isomorphism $\operatorname{Pic}(A) \cong \operatorname{Pic}(A[t])$ iff $A$ is semi-normal. Hence for a non semi-normal ring $A$ there are finite rank projective $A[t]$-modules that are not extended from $A$. Hence by Lemma \ref{BQ}
it follows Theorem 3.11 in \cite{nuss} do not hold for ordinary connections on commutative unital rings. 

Note that Nuss defines a connection using the ideal of the diagonal 
$I \subseteq S\otimes_R S$ and not the module of Kahler differentials $\Omega^1_{S/R}:=I/I^2$. He moreover uses the Amitsur complex to define the curvature of a connection and not the ordinary deRham complex.



\begin{example} Complex line bundles on the projective
  line. \end{example}

Let $\C$ be the field of complex numbers and consider the projective
line $\P^1$ over $\C$. Let $\R$ be the field of real numbers. It
follows the underlying real variety $\P^1(\R)$ of $\P^1$ is isomorphic
to the real $2$-sphere $S^2_\R=\Spec(A)$ where $A=\R[x,y,z]/x^2+y^2+z^2-1$. Any linebundle $\O(d)$ on $\P^1$
gives rise to a locally free $A$-module $P$ of rank $2$. In
\cite{maa2} we proved the linebundle $\O(d)$ on $\P^1$ does not have a
classical connection
\[ \nabla: \O(d)\rightarrow \Omega_{\P^1}^1 \otimes \O(d) \]
for any $d\geq 1$. The underlying real rank $2$ projective $A$-module
$P$ corresponding to  $\O(d)(\R)$ has a classical connection  
\[ \nabla': P\rightarrow \Omega_A^1 \otimes_A P \]
and by
the above results we may give explicit formulas for $\nabla'$ using the
notion of a projective basis and Theorem \ref{basisconnection}.

\section{A formula for the curvature of a connection}

Let in this section $A$ be a commutative unital ring over a fixed
subring $\Z$ and left $E$ be a finitely generated left projective
$A$-module.
Let $B=\{e_1,..,e_n\}$ and $B^*=\{x_1,..,x_n\}$ be a projective basis
for $E$ and let $p:A\{u_1,..,u_n\}\rightarrow E$ be defined by
$p(u_i)=e_i$. It has a left $A$-linear section $s$ defined by
\[ s(e)=\sum_{i=1}^n x_i(e)u_i.\]
We get from Theorem \ref{basisconnection} a connection
\[ \nabla:\Der_\Z(A)\rightarrow \End_\Z(E) \]
defined by
\[ \nabla(\delta)(e)=\sum_{i=1}^n \delta(x_i(e))e_i.\]
The aim of this section is to give a general formula for the trace of
the curvature of the connection $\nabla$ using the notion of a
projective basis. As a consequence we prove
the first Chern class $c_1(E)$ with values in Lie-Rinehart cohomology
and algebraic DeRham cohomology is zero for any finitely generated
projective $A$-module $E$.

\begin{definition}
Let
\[ M=M(p,s)=
\begin{pmatrix} x_1(e_1) & x_1(e_2) & \cdots & x_1(e_n) \\
                x_2(e_1) & x_2(e_2) & \cdots & x_2(e_n) \\
                   \vdots&  \vdots  &  \vdots & \vdots \\
                x_n(e_1) & x_n(e_2) & \cdots & x_n(e_n)
\end{pmatrix}_{B,B^*}
\]
be the \emph{fundamental matrix} of the split surjection
$p$ with respect to the projective basis $B,B^*$.
\end{definition}

Let $\delta\in \Der_\Z(A)$ be a derivation
and let
\[ e=a_1e_1+\cdots +a_ne_n\in E\]
be any element and write
\[e_B= 
\begin{pmatrix} a_1 \\
                a_2 \\
                \vdots \\
                a_n 
\end{pmatrix}_B.
\]

We get the following calculation:
\[ \nabla(\delta)(e)=\]
\[a_1\nabla(\delta)(e_1)+\delta(a_1)e_1+\cdots
+a_n\nabla(\delta)(e_n)+\delta(a_n)e_n.\]
By definition
\[ \nabla(\delta)(e_k)=\delta(x_1(e_k))e_1+\delta(x_2(e_k))e_2+\cdots
+\delta(x_n(e_k))e_n.\]
We write this in ``matrix form'' as follows:
\[
\nabla(\delta)(e_k)=
\begin{pmatrix} \delta(x_1(e_k)) \\
                \delta(x_2(e_k)) \\
                  \vdots \\
                \delta(x_n(e_k))
\end{pmatrix}_{B }.
\]

We will write down explicit formulas for the connection $\nabla$ in
matrix-notation using the above notation. 
By definition
\[\nabla(\delta)(e_k)= \delta(x_1(e_k))e_1+\cdots
+\delta(x_n(e_k))e_n=\]
\[
\begin{pmatrix} \delta(x_1(e_k)) \\
                \delta(x_2(e_k)) \\
              \cdots              \\
               \delta(x_n(e_k)) 
\end{pmatrix}_{B }
=[ \delta(x(e_k))]_{B }.
\]
We get
\[ \nabla(\delta)(e)=\]
\[a_1[\delta(x(e_1))]_{B }+\delta(a_1)e_1+\cdots
+a_n[\delta(x(e_n))]_{B }+\delta(a_n)e_n=\]
\[ ([D_\delta]_B +[\delta(M)]_B)e_B \]
where we use the following notation:
\[ [D_\delta]_B(e_B)=
\begin{pmatrix} \delta(a_1) \\
               \delta(a_2) \\
 \vdots \\
\delta(a_n)
\end{pmatrix}_{B }
\]
and
\[ [\delta(M)]_B=
\begin{pmatrix} \delta(x_1(e_1)) & \delta(x_1(e_2)) & \cdots & \delta(x_1(e_n)) \\
                \delta(x_2(e_1)) & \delta(x_2(e_2)) & \cdots & \delta(x_2(e_n)) \\
                   \vdots&  \vdots  &  \vdots & \vdots \\
                \delta(x_n(e_1)) & \delta(x_n(e_2)) & \cdots & \delta(x_n(e_n))
\end{pmatrix}_{B }.
\]
It follows we have described the endomorphism $\nabla(\delta)$
completely in terms of the fundamental matrix $M=M(p,s)$:
\[ \nabla(\delta)=[D_{\delta}]_B+[\delta(M)]_B\]
and the projective basis $\{B,B^*\}$.
We want to calculate the curvature $R_\nabla$ using $B$ and
$B^*$. This is a long but straight forward calculation which we now
present. We get
\[ \nabla(\eta)\nabla(\delta)(e)=\]
\[\nabla(\eta)(a_1\sum_{i=1}^n \delta(x_i(e_1))e_i +a_2\sum_{i=1}^n
\delta(x_i(e_2))e_i+\cdots +a_n\sum_{i=1}^n \delta(x_i(e_n))e_i) \]
\[ +\delta(a_1)e_1+\delta(a_2)e_2+\cdots+\delta(a_n)e_n)=\]
\[ \sum_{i=1}^na_1\delta(x_i(e_1))\nabla(\eta)(e_i) +
\eta(a_1\delta(x_i(e_1)))e_i \]
\[+ \sum_{i=1}^na_2\delta(x_i(e_1))\nabla(\eta)(e_i) +
\eta(a_2\delta(x_i(e_1)))e_i+\cdots + \]
\[ \sum_{i=1}^na_n\delta(x_i(e_1))\nabla(\eta)(e_i) +
\eta(a_n\delta(x_i(e_1)))e_i \]
\[ +\sum_{i=1}^n \delta(a_i)\nabla(\eta)(e_i)+\eta(\delta(a_i))e_i=\]
\[ \nabla(\eta)(e_1)(\sum_{i=1}^n a_i\delta(x_1(e_i)))+ \]
\[ \nabla(\eta)(e_2)(\sum_{i=1}^n a_i\delta(x_2(e_i)))+\cdots + \]
\[ \nabla(\eta)(e_n)(\sum_{i=1}^n a_i\delta(x_n(e_i)))+ \]
\[ e_1(\sum_{i=1}^n \eta(a_i\delta(x_1(e_i))))+e_2(\sum_{i=1}^n
\eta(a_i\delta(x_2(e_i)))+ \cdots +\]
\[ e_n(\sum_{i=1}^n \eta(a_i\delta(x_n(e_i)))+\]
\[ \sum_{i=1}^n\delta(a_i)\nabla(\eta)e_i+\eta(\delta(a_i))e_i.\]

When we express the above calculation in matrix-notation we get
\begin{align}
&\label{matrix}\nabla(\eta)\nabla(\delta)(e)= 
\end{align}
\[ [\eta(M)]_B[\delta(M)]_B(e_B)+[\delta(M)]_B[D_\eta]_B(e_B)+\]
\[ [(\eta\circ \delta)(M)]_B(e_B)+[\eta(M)]_B[D_\delta]_B(e_B)+[D_{\eta\circ \delta}]_B(e_B).\]

\begin{theorem} \label{curvature} The following holds:
\[ R_{\nabla}(\delta,\eta)=[[\delta(M)]_B,[\eta(M)]_B] .\]
\end{theorem}
\begin{proof} By definition 
\[ R_\nabla(\delta,\eta)=[\nabla(\delta),\nabla(\eta)]-\nabla([\delta,
\eta])=\]
\[ \nabla(\delta)\circ
\nabla(\eta)-\nabla(\eta)\circ\nabla(\delta)-\nabla([\delta,\eta]).\]
Since $\nabla(\delta)=[D_{\delta}]_B+[\delta(M)]_B$ we get using
Formula \ref{matrix}  the following calculation:
\[ R_\nabla(\delta,\eta)(e)=\]
\[
\nabla(\delta)\nabla(\eta)(e)-\nabla(\eta)\nabla(\delta)(e)-\nabla([\delta,\eta])(e)=\]
\[ [\delta(M)]_B[\eta(M)]_B(e_B)+[\eta(M)]_B[D_\delta]_B(e_B) \]
\[ +[(\delta\circ \eta)(M)]_B(e_B)+[\delta(M)]_B[D_\eta]_B(e_B)+[D_{\delta\circ
  \eta}]_B(e_B)\]
\[-[\eta(M)]_B[\delta(M)]_B(e_B)-[\delta(M)]_B[D_\eta]_B(e_B)-[(\eta\circ
\delta)(M)]_B(e_B) \]
\[ -[\eta(M)]_B[D_\delta]_B(e_B)-[D_{\eta\circ \delta}]_B(e_B)\]
\[ -[D_{[\delta,\eta]}]_B(e_B)-[ [\delta,\eta](M)]_B(e_B)=\]
\[ [ [\delta(M)]_B,[\eta(M)]_B](e_B)\]
and the claim of the Theorem follows.
\end{proof}

Since the curvature matrix $R_\nabla(\delta\wedge \eta)$ is the
commutator of two matrices it follows
\[ det(R_\nabla(\delta\wedge \eta))=0 \]
when $n$ is an odd number.

Note: The trace of the commutator $R=[[\delta(M)]_B, [\eta(M)]_B]$ as an element of $\End_A(E)$ is non zero in general.
If $e_1,\ldots ,e_n, x_1,\ldots ,x_n $ is a projective basis the trace of $R$ is given as follows:
\[ tr(R)=\sum_{i=1}^n x_i(R(e_i)).\]
The commutator $[[\delta(M)]_B, [\eta(M)]_B]$ has trace zero as an element of $\End_A(A^n)$.
See \cite{maa0} for an explicit example of a connection $\nabla$ where $tr(R_\nabla)\neq 0$.

\section{Examples: Algebraic connections on ellipsoids}

In this section we use the notions introduced in the previous section
to give an explicit construction of connections on modules of Kahler-differentials on ellipsoids.

Let in the following $A$ be any $\mathbf{Q}$-algebra and let
$B=A[x_1,..,x_k]$ be the polynomial ring in $k$ variables over $A$. Let
\[ H=x_1^{p_1}+\cdots +x_k^{p_k}-1\in B\]
and consider the $A$-algebra $C=B/H$. In this section we apply the
construction in the previous section to give explicit formulas for
connections on $\Omega=\Omega_{C/A}$. Let $d:C\rightarrow \Omega$ be
the universal derivation. Let
\[ dH=p_1x_1^{p_1-1}dx_1+\cdots +p_kx_k^{p_k-1}dx_k.\]
It follows there is an isomorphism
\[ \Omega=C\{dx_1,..,dx_k\}/dH \]
of left $C$-modules. Let $F=C\{e_1,..,e_n\}$ be the free $C$-module on
the basis $\{e_1,..,e_n\}$. Let 
\[ G=p_1x_1^{p_1-1}e_1+\cdots +p_kx_k^{p_k-1}e_k \in F\]
and let $Q$ be the left $C$-module spanned by $G$.
We get an exact sequence of left
$C$-modules
\[ 0\rightarrow Q\rightarrow F \rightarrow^p \Omega \rightarrow 0\]
where the map $p$ is defined as follows:
\[ p(e_i)=\overline{dx_i}.\]
Since the module $\Omega$ is projective there is a $C$-linear
splitting $s$ of $p$. Define the following map $s:\Omega \rightarrow
F$:
\[ s(\overline{dx_i})=e_i-\frac{1}{p_i}x_i G.\]

\begin{lemma} \label{ellipsoidsection} The map $s$ is a left $C$-linear section of $p$.
\end{lemma}
\begin{proof} We prove the map $s$ is well defined: 
\[ s(dH)=s(p_1x_1^{p_1-1}\overline{dx_1}+\cdots
+p_kx_k^{p_k-1}\overline{dx_k})=\]
\[p_1x_1^{p_1-1}(e_1-\frac{1}{p_1}x_1G)+\cdots
+p_kx_k^{p_k-1}(e_k-\frac{1}{p_k}x_kG)=\]
\[p_1x_1^{p_1-1}e_1+\cdots +p_kx_k^{p_k-1}e_k-x_1^{p_1}G-\cdots
-x_k^{p_k}G=\]
\[ G-(x_1^{p_1}+\cdots +x_k^{p_k})G=G-G=0.\]
It follows
\[ s(dH)=0 \]
hence the map $s$ is well defined. By definition 
\[ s(\omega)=s(a_1\overline{dx_1}+\cdots +a_k\overline{dx_k})=\]
\[ a_1(e_1-\frac{1}{p_1}x_1G)+\cdots + a_k(e_k-\frac{1}{p_k}x_kG)=\]
\[ a_1e_1+\cdots +a_ke_k + yG.\]
it follows
\[ p(s(\omega))=p(a_1e_1+\cdots +a_ke_k+yG)= \]
\[ p(a_1e_1+\cdots +a_ke_k) + p(yG)= \omega.\]
The Lemma is proved.
\end{proof}

Let $e_i^*:F\rightarrow A$ be coordinate functions on $F$ and put
$x_i=e_i^*\circ s$. Let $w_i=\overline{dx_i}$ for $i=1,..,k$. Let
$B=\{w_1,..,w_k\}$ and $B^*=\{x_1,..,x_k\}$. It follows $\{B,B^*\}$ is
a projective basis for the projective module $\Omega$. 
Define the following map:
\[ \nabla: \Omega \rightarrow \Omega\otimes_A \Omega \]
by
\[ \nabla(w)=\sum_i d(x_i(w))\otimes w_i .\]

\begin{theorem} \label{omega}  The map $\nabla$ is a connection on $\Omega$.
\end{theorem}
\begin{proof} By Lemma \ref{ellipsoidsection} it follows the sets
  $\{B,B^*\}$ form a projective basis for $\Omega$. It follows from
  Theorem \ref{basisconnection} the map $\nabla$ is a connection on $\Omega$.
\end{proof}

In the following we calculate explicitly connections and their
curvature on a class of ellipsoids.

\begin{example} Connections on the two-sphere.\end{example}

Let in this example $K$ be a field of characteristic different from
two. Let $f=x_1^2+x_2^2+x_3^2-1 \in K[x_1,x_2,x_3]$.
Let $A=K[x_1,x_2,x_3]/f$. It follows $S=\Spec(A)$ is the two-sphere over $K$.
Let $\Omega=\Omega_{A/K}$ be the module of Kahler differentials of $A$
relative to $K$. It follows $\Omega=A\{dx_1,dx_2,dx_3\}/H$ where
$H=x_1dx_1+x_2dx_1+x_3dx_3$. Let $G=x_1u_1+x_2u_2+x_3u_3$ where
$A\{u_1,u_2,u_3\}$ is the free $A$-module of rank $3$. Since $\Der_K(A)$
is a non-trivial rank two locally free $A$-module it follows $\Omega$
is a non-trivial rank two locally free $A$-module.

We get an exact
sequence
\[ 0\rightarrow (G)\rightarrow A\{u_1,u_2,u_3\}\rightarrow^p \Omega
\rightarrow 0 \]
where $p(u_i)=\overline{dx_i}$. The section $s$ defined in Lemma
\ref{ellipsoidsection} defines a projective basis $B,B^*$ and a
fundamental matrix $M=M(p,s)$ for $\Omega$ with respect to the split
surjection $p$.

Let $T=\Der_K(A)$. It follows $T$ is
generated by the following derivations:
\[ \PARone=x_2\parone-x_1\partwo \]
\[ \PARtwo=x_3\parone-x_1\parthree \]
and
\[ \PARthree=x_3\partwo-x_2\parthree.\]
let $\omega=a_1dx_1+a_2dx_2+a_3dx_3\in \Omega$ and let $\partial
\in \End_K(A)$. We define
\[ D_{\partial}(\omega)=\partial(a_1)dx_1+\partial(a_2)dx_2+\partial(a_3)dx_3\in \Omega.\]
Define the following elements:
\[
\nabla(\PARone)= [D_{\PARone}]_B+
\begin{pmatrix} -2x_1x_2      & x_1^2-x_2^2 & -x_2x_3 \\
                x_1^2-x_2^2   &  2x_1x_2    & x_1x_3  \\
                -x_2x_3       &  x_1x_3     &  0 
\end{pmatrix}
=[D_{\PARone}]_B+[\PARone(M)]_B
\]

\[
\nabla(\PARtwo)= [D_{\PARtwo}]_B+
\begin{pmatrix} -2x_1x_3      & -x_2x_3 & x_1^2-x_3^2 \\
                 -x_2x_3   &    0       & x_1x_2      \\
                x_1^2-x_3^2        &  x_1x_2      &  2x_1x_3
\end{pmatrix}
=[D_{\PARtwo}]_B+[\PARtwo(M)]_B
\]

and

\[
\nabla(\PARthree)= [D_{\PARthree}]_B+
\begin{pmatrix}    0      & -x_1x_3    & x_1x_2       \\
                -x_1x_2   &  -2x_2x_3       & x_2^2-x_3^2  \\
                x_1x_2       &  x_2^2-x_3^2     &  2x_2x_3
\end{pmatrix}
=[D_{\PARthree}]_B+[\PARthree(M)]_B
\]

\begin{corollary} The maps $\nabla(\PARone). \nabla(\PARtwo)$ and
  $\nabla(\PARthree)$ define a connection
\[ \nabla: \Der_K(A)\rightarrow \End_K(\Omega) .\]
\end{corollary}
\begin{proof} One checks the given formulas are the formulas one gets
  when one makes the connection in Theorem \ref{omega} explicit.
\end{proof}

The curvature of the connection $\nabla$ is the following map:
\[ R_\nabla: \Der_K(A)\wedge \Der_K(A)\rightarrow \End_K(\Omega) \]
\[ R_\nabla(\delta\wedge
\eta)=[\nabla(\delta),\nabla(\eta)]-\nabla([\delta,\eta]).\]

We calculate $R_\nabla(\partial_i\wedge \partial_j)$ using
Theorem \ref{curvature}:
The fundamental matrix $M=M(p,s)$ is the following matrix:
\[
M=
\begin{pmatrix} 1-x_1^2 & -x_1x_2 & -x_1x_3 \\
               -x_1x_2 & 1-x_2^2 & -x_2x_3 \\
                -x_1x_3 & -x_2x_3 & 1-x_3^2
\end{pmatrix}.
\]
It follows
\[ [\partial_1(M)]_B=
\begin{pmatrix} -2x_1x_2 & x_1^2-x_2^2 & -x_1x_3 \\
             x_1^2-x_2^2 & 2x_1x_2 & x_1x_3 \\
           -x_2x_3 & x_1x_3 & 0
\end{pmatrix}
\]
and
\[ [\partial_2(M)]_B=
\begin{pmatrix} -2x_1x_3  & -x_2x_3 & x_1^2-x_3^2 \\
                 -x_2x_3  & 0       & x_1x_2 \\
             x_1^2-x_3^2 & x_1x_2  & 2x_1x_3 
\end{pmatrix}
\]
One calculates using Theorem \ref{curvature}
\[ R_\nabla(\partial_1\wedge \partial_2)=[[\partial_1(M)]_B,[\partial_2(M)]_B]=
\begin{pmatrix} 0 & x_1x_3 & -x_1x_2 \\
            -x_1x_3 & 0 & x_1^2 \\
        x_1x_2 & -x_1^2 & 0 
\end{pmatrix}.
\]

Let 
\[A=
\begin{pmatrix} a_{11} & a_{12} & a_{13} \\
a_{21}  & a_{22} & a_{23} \\
a_{31} & a_{32}  & a_{33}
\end{pmatrix}
\]
be any matrix with $a_{ij}$ independent variables.
Consider the characteristic polynomial $P_A(\lambda)=det(\lambda I-A)$ 
where $I$ is the $3 \times 3$ identity matrix.
It follows
\[ P(\lambda)=\lambda^3-tr(A) \lambda^2 +p_A \lambda -det(A) \]
where
\[p_A=a_{11}a_{22}+a_{11}a_{33}+a_{22}a_{33}+ a_{21}a_{12}+
a_{31}a_{13}+a_{32}a_{23} .\]
Let $A=R_\nabla(\partial_1\wedge \partial_2)$.
It follows $p_A=-x_1^2\neq 0$.

By \cite{maa1} the connection $\nabla$ define a Chern-class
\[ c_1(\Omega)\in \H^2(\Der_K(A),A) \]
where $\H^2(\Der_K(A),A)$ is the second Lie-Rinehart cohomology of $\Der_K(A)$.

\begin{corollary}\label{flat} The following holds:
\begin{align}
&\label{c1}\text{$\nabla$ is non-flat.}\\
&\label{c2}\text{$tr(R_{\nabla})=0$.}\\
&\label{c3}\text{$c_1(\Omega)=0$ in $\H^2(\Der_K(A),A)$.}
\end{align}
\end{corollary}
\begin{proof} Claim \ref{c1} follows since the map
  $R_{\nabla}(\partial_1 \wedge \partial_2)$ is a non-zero element of
  $\End_A(\Omega)$ as one easily checks.
Claim \ref{c2} and  Claim \ref{c3} follows from an explicit calculation. The Corollary is proved.
\end{proof}

Note: We get examples of non-flat algebraic  connections
\[ \nabla:\Der_K(A)\rightarrow \End_K(\Omega) \]
defined over any field $K$ of characteristic different from two.

\begin{example} Non flat algebraic connections on projective modules on the real 2-sphere. \end{example}

Note: By \cite{milnor} the following holds: Asume $k$ is the field of real numbers and assume $I:=(f_1,..,f_l) \subseteq k[x_1,..,x_n]$ is an ideal whose zero set
$X:=Z(I) \subseteq k^n$ define a real smooth simply connected manifold $X(k)$. Let $A:=k[x_i]/I$ and let $E$ be a finite rank projective $A$-module whose coorresponding real smooth vector bundle
$E(k)$ on $X(k)$ is non trivial. It follows by \cite{milnor}, page 294 that $E(k)$ has no flat smooth connection. A proof of this claim may be found in \cite{taubes}. Hence if $\nabla$ is a flat algebraic connection on $E$, 
it follows $\nabla$ gives rise to a smooth connection $\nabla_k$ on $E(k)$ since $\nabla$ is defined using polynomial functions. From \cite{milnor} it follows $E(k)$ is trivial, a contradiction.
Hence if $E(k)$ is non  trivial it follows $E$ has no flat algebraic connection. The tangent bundle $T_{S^2}$ and cotangent bundle $\Omega_{S^2}$ on the real 2-sphere $S^2$ 
are non-trivial as smooth vector bundles, hence by \cite{milnor} and \cite{taubes} it follows $T_{S^2}$ and $\Omega_{S^2}$ have no flat algebraic connections.  Hence Corollary \ref{flat} gives an example
of a non flat algebraic connection 

\[ \nabla: T_{S^2} \rightarrow \End_k(\Omega_{S^2}) \]

on the cotangent bundle $\Omega_{S^2}$. The contangent bundle $\Omega_{S^2}$ is non-trivial but its Chern classes are trivial since $\Omega_{S^2}$ is stably trivial. This follows from the Whitney sum formula.
In the paper \cite{maa01} I prove that the connection $\nabla$ from Corollary \ref{flat} gives rise to a non-trivial characteristic class

\[ nc_1(\Omega_{S^2}) \in \Ext^1(T_{S^2}, \End_A(\Omega_{S^2})).\]

The class $nc_1(\Omega_{S^2})$ is trivial iff $\Omega_{S^2}$ has a flat algebraic connection.

\begin{example} Complex projective manifolds and holomorphic
  connections. \end{example}

In \cite{atiyah} the following is conjectured: Let $X$ be a complex
projective manifold and let $\E$ be a holomorphic bundle with a
holomorphic connection $\nabla$. Then there exist a holomorphic flat
connection $\nabla'$. In the affine case as indicated above many
``naturally occuring'' algebraic connections are non-flat.

\section{Higher order differential operators}

In the litterature one finds many papers devoted to the construction
of explicit formulas for connections on 
maximal Cohen-Macaulay modules on isolated hypersurface singularities
(See \cite{maa1} for one approach to this problem using the
Kodaira-Spencer map and the Atiyah class). In this section we use the
notion of a projective basis on a finitely generated projective
$A$-module to give explicit formulas for $l$-connections
\[ \nabla^l:E\rightarrow P^l\otimes_A E \]
where $P^l$ is the \emph{$l$'th module of principal parts} of the ring
$A$. The connections $\nabla^l$ are non-flat in general. We also
consider the notion of a stratification and give explicit examples of
projective finitely generated module where the canonical connection
induced by a projective basis does not give rise to a stratification. 
The obstruction to this is given by the $(l,k)$-curvature $K^{(l,k)}$
of the $l$-connection $\nabla^l$.

Let in the following $\Z$ be a fixed commutative unital base ring and
let $A$ be a commutative $\Z$-algebra. Let $E$ be a left $A$-module.

The existence of a flat connection
\[ \nabla: \Der_\Z(A)\rightarrow \End_\Z(E) \]
on $E$ induce a ring homomorphism
\[ \rho: \diff_\Z(A)\rightarrow\End_\Z(E) \]
where $\diff_\Z(A)$ is the \emph{small ring of differential operators}
of $A$. The module $E$ is by definition a left $A$-module and the induced
structure $\rho$ is a lifting of the left $A$-module structure to a
left $\diff_\Z(A)$-module structure. The obstruction to this lifting is the
curvature of the connection $\nabla$.
The ring $\diff_\Z(A)$ is the associative subring of $\End_\Z(A)$ generated by
$\Der_\Z(A)$. In general the ring $\diff_\Z(A)$ is a strict subring of
the ring $\Diff_\Z(A)$ - the \emph{ring of differential operators}
of $A$. In the case where $A$ is a regular $k$-algebra of finite type
over $k$ where $k$ is a field of characteristic zero one has an
equality
\[ \diff_\Z(A)=\Diff_\Z(A) .\]
In the case when $A$ is non-regular there is a strict inclusion of
rings $\diff_\Z(A)\subseteq \Diff_\Z(A)$. The higher order connections $\nabla^l$
we construct are related to the following problem:
One wants to lift the left $A$-module structure on $E$ to a left
$\Diff_\Z(A)$-module structure
\[ \rho:\Diff_\Z(A) \rightarrow \End_\Z(E) \]
on $E$. The obstruction to such a lifting is given by the notion of \emph{generalized curvature} of the
connection $\nabla^l$.

The ring $\Diff_\Z(A)$ has a filtration $\Diff^l_\Z(A)$ by
degree of differential operators and there is an isomorphism
\[ \Diff^l_\Z(A)\cong \Hom_A(P^l_{A/\Z},A) \]
where $P^l=P^l_{A/\Z}=A\otimes_\Z A/I^{l+1}$ is the $l$'th module of principal parts of
$A$. Here $I\subseteq A\otimes_\Z A$ is the kernel of the canonical
multiplication map. 
Define the following map
\[ \partial^l:A\rightarrow P^l\]
by
\[ \partial^l(a)=\overline{1\otimes a} .\]
There is a canonical projection map
\[ p_l:P^l \rightarrow P^{l-1} \]
and an equality $p_l \circ \pl =\partial^{l-1}$ of maps.
Let $E,F$ be left $A$-modules and define the following:
\[ \Diff^{-1}_\Z(E, F)=0\]
and
\[ \Diff^l_\Z(E,F)=\{ \partial \in \Hom_\Z(E,F): [\partial, a]\in
\Diff^{l-1}_\Z(E,F)\text{ for all $a\in A$} \} .\]
By definition $\Diff^l_\Z(A)=\Diff^l_\Z(A,A)$ and 
\[ \Diff_\Z(A)=\cup_{l\geq -1}\Diff^l_\Z(A).\]
We get a filtration of left and right $A$-modules
\[ 0=\Diff^{-1}_\Z(E,F)\subseteq \Diff^0_\Z(E,F)\subseteq \cdots \subseteq
\Diff^l_\Z(E,F) \subset \cdots \subseteq \Diff_\Z(E,F) \]
where
\[ \Diff_\Z(E,F) =\cup_{l\geq -1}\Diff^l_\Z(E,F).\]

Consider the following map
\[ d:A\rightarrow A\otimes_\Z A\]
defined by
\[ d(a)=1\otimes a-a\otimes 1.\]

Let $I_l=\{1,2, \cdots ,l\}$.
Given $a\in A$, let $\phi_a \in \End_\Z(A)$ be the endomorphism which
is multiplication with $a$. Let $\partial \in \End_\Z(A)$.

\begin{lemma} \label{diffformulas} The following holds:
\begin{align}
&\label{m1} d(a_1)\cdots d(a_l)=\sum_{H\subseteq I_l}(-1)^{card(H)}(\prod_{i\in
  H}a_i)\otimes (\prod_{i\notin H}a_i)\in A\otimes_\Z A. \\
&\label{m2}[\cdots [\partial, \phi_{a_1}]\cdots
]\phi_{a_1}](a)=\sum_{H\subseteq I_l}(-1)^{card(H)}(\prod_{i\in
  H}a_i)\partial(\prod_{i\notin H}a_i)a)
\end{align}
\end{lemma}
\begin{proof} The proof of the Lemma follows from \cite{groth}, Proposition 16.8.8 using induction.
\end{proof}

Let $\partial \in \Hom_\Z(E,E)$ be a $\Z$-linear endomorphism and define
the following map:
\[ \phi_{\partial}:A\otimes_\Z A\otimes_A E\rightarrow E\]
by
\[ \phi_{\partial}(a\otimes b\otimes e)=a\partial(be).\]
Let $I\subseteq A\otimes_\Z A$ be the kernel of the multiplication map.

\begin{proposition} The following holds:
\[ \partial \in \Diff^l_\Z(E)\text{ if and only if }
\phi_{\partial}(I^{l+1}(A\otimes_\Z A\otimes_A E)=0.\] 
\end{proposition}
\begin{proof}
Assume $\partial \in \Diff^l_\Z(A)$. It follows from Lemma
\ref{diffformulas}
\[ \phi_{\partial}(\prod_{i=1}^{l+1}(1\otimes a_i-a_i\otimes 1)\otimes
e)=\]
\[phi_{\partial}(\sum_{H\subseteq I_{l+1}}(-1)^{card(H)}(\prod_{i\in
  H}a_i)\otimes (\prod_{i\notin H}a_i)\otimes e)=\]
\[= \sum_{H\subseteq I_{l+1}}(-1)^{card(H)}(\prod_{i\in
  H})\partial((\prod_{i\notin H}a_i)e)=\]
\[ [\cdots [\partial ,\phi_{a_1}]\cdots ]\phi_{a_{l+1}}](e)=0\]
since $\Diff^{-1}_\Z(A)=0$. It follows
$\phi_{\partial}(I^{l+1}(A\otimes A\otimes E))=0$. The converse
statement is proved in a similar way and the Proposition is proved.
\end{proof}

\begin{lemma} \label{diff} The following holds:
\[ \pl\in \Diff^l_\Z(A,P^l).\]
\end{lemma}
\begin{proof} Define the following map
\[ \partial:A\rightarrow A\otimes_\Z A\]
by
\[ \partial(a)=1\otimes a.\]
One gets the formula
\[ [\cdots [\partial, \phi_{a_1}]\cdots ]\phi_{a_{l+1}}](x)=d(a_1)\cdots
d(a_{l+1})(1\otimes x) \]
where $d(a)=1\otimes a-a\otimes 1$ and where the product is in
$A\otimes_\Z A$. 
There is a natural projection map
\[ p_l:A\otimes_\Z A\rightarrow P^l \]
deifned by
\[ p_l(a\otimes b)=\overline{a\otimes b}.\]
We get $\pl =p_l\circ \partial$. 
It follows
\[ [\cdots [\pl, \phi_{a_1}]\cdots ]\phi_{a_{l+1}}](x)= \]
\[ \overline{[\cdots [\partial , \phi_{a_1}]\cdots ]\phi_{a_{l+1}}](x)}= \]
\[ \overline{ d(a_1)d(a_2)\cdots d(a_{l+1})(1\otimes x)}=0 \]
hence
\[ [\cdots [\pl,\phi_{a_1}]\cdots ]\phi_{a_{l+1}}] \in \Diff^{-1}_\Z(A,P^l)=0 .\]
It follows
\[ \pl \in \Diff^l_\Z(A,P^l) \]
and the Lemma is proved.
\end{proof}

Let $E$ be a left projective $A$-module with projective basis
$B=\{e_1,..,e_n\}$ and $B^*=\{x_1,..,x_n\}$ and consider the
connection
\[ \nabla:E\rightarrow I \otimes_A E \]
defined by
\[ \nabla(e)=\sum_{i=1}^n d(x_i(e))\otimes e_i \]
where 
\[ d:A\rightarrow I \]
is a derivation: $ d\in \Der_\Z(A,I)$. By Theorem \ref{basisconnection} it
follows $\nabla$ is a connection on $E$.

\begin{lemma} The following holds:
\[ \nabla \in \Diff^1_\Z(E, I\otimes_A E) .\]
\end{lemma}
\begin{proof} Let $a\in A$ and consider the following map:
\[ [\nabla, a]:E\rightarrow \Omega\otimes_A E.\]
It follows
\[ [\nabla,a](e)=\nabla(ae)-a\nabla(e)=a\nabla(e)+d(a)\otimes e
-a\nabla(e)=\]
\[ d(a)\otimes e.\]
We get
\[ [\nabla,a](be)=d(a)\otimes be=b(d(a)\otimes e)=b[\nabla, a](e).\]
It follows
\[ [\nabla, a]\in \Hom_A(E, I\otimes
E)=\Diff^0_\Z(E, I \otimes_A E) \]
and hence
\[ \nabla\in \Diff^1_\Z(E, I \otimes_A E).\]
The Lemma is proved.
\end{proof}

Define the following map
\[ \nabla^l:E\rightarrow P^l\otimes_A E\]
by
\[ \nabla^l(e)=\sum_{i=1}^n \pl(x_i(e))\otimes e_i.\]
Since $\pl\in \Diff_\Z^l(A,P^l)$ it follows
\[ \nabla^l\in \Hom_\Z(E, P^l\otimes_A E).\]

\begin{theorem} \label{higher} The following holds:
\begin{align}
&\label{h1}  \nabla^l\in \Diff^l_\Z(E,P^l\otimes_A E).\\
&\label{h2} (p_l\otimes 1)\circ \nabla^l=\nabla^{l-1}.
\end{align}
\end{theorem}
\begin{proof}
By Lemma \ref{diff} the following formula holds:
\[ [\cdots [\nabla^l,\phi_{a_1}]\cdots ]\phi_{a_{l+1}}](e)=\sum_{i=1}^n[\cdots
[\pl,\phi_{a_1}]\cdots ]\phi_{a_{l+1}}](x_i(e))\otimes e_i=\]
\[ \sum_{i=1}^n \overline{d(a_1)\cdots d(a_{l+1})(x_i(e))}\otimes e_i=0\]
It follows 
\[ [\cdots [\nabla^l,\phi_{a_1}]\cdots ]\phi_{a_{l+1}}]\in \Diff^{-1}_\Z(E,P^l\otimes_A
E)=0  \]
hence
\[ \nabla^l\in \Diff^l_\Z(E, P^l\otimes_A E) \]
and the Proposition is proved.
\end{proof}

We get for any left projective $A$-module $E$ and all $l\geq 1$ commutative diagrams
\[ 
\diagram E \rto^{\nabla^l} \drto_{\nabla^{l-1}} & P^l\otimes_A E \dto^{p_l\otimes 1} \\
   & P^{l-1}\otimes_A E 
\enddiagram \]
of differential operators.

\begin{definition} The map
\[ \nabla^l:E\rightarrow P^l\otimes_A E\]
is the \emph{$l$-connection} associated to the projective
basis $B,B^*$. 
\end{definition}

Note: There $l$'th module of principal parts $P^l=A\otimes_\Z
A/I^{l+1}$ is a commutative unital ring in an obvious way and there
is a multiplicative unit $1\in P^l$. One may define the following map
\[ \rho^l :E\rightarrow P^l\otimes_A E \]
by
\[ \rho^l(e)=1\otimes e.\]
The map $\rho^l$ is called the \emph{universal differential operator}
for $E$ of order $l$. One checks $\rho^l\in \Diff^l_\Z(E,P^l\otimes_A
E)$.  The map $\rho^l$ induce an isomorphism
\begin{align}
&\label{repr} \Diff^l_\Z(E,F)\cong \Hom_A(P^l\otimes_A E, F) 
\end{align}
of left and right $A$-modules.

Given an $l$'th order differential operator $\partial \in
\Diff^l_\Z(A)$ one gets by Formula \ref{repr} an $A$-linear map
$\phi_{\partial}:P^l \rightarrow A$
We get using $\nabla^l$ a map
\[ \rho(\partial):E\rightarrow E \]
defined by
\[ \rho(\partial)(e)=\phi_{\partial}\otimes
1(\nabla^l(e))=\sum_{i=1}^n\phi_{\partial}(\pl(x_i(e)))e_i=\sum_{i=1}^n \partial(x_i(e))e_i.\]
This defines a map of left $A$-modules
\[ \rho^l:\Diff^l_\Z(A)\rightarrow \Diff^l_\Z(E) \]
for all$l\geq 1$. By Theorem \ref{higher} we get commutative diagrams
of maps
\[ 
\diagram  \Diff^l_\Z(A)\dto^j  \rto^{\rho^l} & \Diff^l_\Z(E) \dto^i  \\
          \Diff^{l+1}(A) \rto^{\rho^{l+1}} & \Diff^{l+1}_\Z(E) 
\enddiagram \]
for all $l\geq 1$ with $i,j$ the canonical inclusion maps. This induce
a canonical map
\begin{align}
&\label{canonical} \rho:\Diff_\Z(A)\rightarrow \Diff_\Z(E) 
\end{align}

\begin{example} \label{TWOsphere} Connections on the two-sphere.\end{example}

Let in this example $K$ be a field of characteristic zero. 
Let $f=x_1^2+x_2^2+x_3^2-1$ be in $K[x_1,x_2,x_3]$. Let $A=K[x_1,x_2,x_3]/f$. It follows
$S=\Spec(A)$ is the two-sphere over $K$.
Let $\Omega=\Omega_{A/K}$ be the module of Kahler differentials of $A$
relative to $K$. It follows $\Omega=A\{dx_1,dx_2,dx_3\}/H$ where
$H=x_1dx_1+x_2dx_1+x_3dx_3$. Let $G=x_1u_1+x_2u_2+x_3u_3$ where
$A\{u_1,u_2,u_3\}$ is the free $A$-module of rank $3$. We get an exact
squence
\[ 0\rightarrow (G)\rightarrow A\{u_1,u_2,u_3\}\rightarrow^p \Omega
\rightarrow 0 \]
where $p(u_i)=\overline{dx_i}$. Let $s(dx_i)=u_i-x_iG$. It follows $s$
is a left $A$-linear section of $p$. It is well known $\Omega$ is a
non-free locally free rank two $A$-module.
Let $\Der_K(A)$ be the module of $K$-linear derivations of $A$. It follows $\Der_K(A)$ is
generated by the following derivations:
\[ \PARone=x_2\parone-x_1\partwo \]
\[ \PARtwo=x_3\parone-x_1\parthree \]
and
\[ \PARthree=x_3\partwo-x_2\parthree.\]
Since $S$ is smooth over $K$ it follows $\Diff_K(A)$ is generated as
an associative ring by $\PARone, \PARtwo$ and $\PARthree$.
There is a map $\rho$ induced by the projective basis $B,B^*$ for
$\Omega$ defined as follows:
\[ \rho:\Diff_K(A)\rightarrow \Diff_K(\Omega) \]
defined by
\[
\rho(\partial)(e)=\partial(x_1(e))dx_1+\partial(x_2(e))dx_s+\partial(x_3(e))dx_3.\]
The map $\rho$ is not a morphism of associative rings for the
following reason:
One calculates the following:
\[ \PARone
\circ \PARtwo=x_2x_3\parone^2+x_2\parthree+x_1x_2\parone\circ \parthree
-x_1x_3\parone\circ \partwo +x_1^2\partwo \circ \parthree.\]
One calculates the following:
\[ \rho(\PARone\circ \PARtwo)(dx_1)=\]
\[-2x_1x_3dx_1+x_1x_3dx_2-2x_1x_2dx_3.\]
One calculates the following:
\[ \rho(\PARone)\circ
\rho(\PARtwo)(dx_1)=\nabla(\PARone)(\nabla(\PARtwo)(dx_1))=\]
\[
(x_1^2x_2x_3+3x_2x_3)dx_1+(3x_1x_2^2x_3-x_1x_3)dx_2+(x_1x_2x_3^2+2x_1x_2)dx_3.\]
Hence
\[ \rho(\PARone \circ \PARtwo)\neq \rho(\PARone)\circ \rho(\PARtwo).\]
It follows $\rho$ is not a morphism of associative rings.

From  Corollary \ref{flat} it also follows the map $\rho$ is not a ring
homomorphism:  Since the connection
\[ \nabla:\Der_K(A)\rightarrow \End_K(\Omega) \]
is non-flat the $A$-module structure on $\Omega$ does not lift to
a $\diff_B(A)$-module structure on $E$ hence it does not lift to a
$\Diff_B(A)$-module structure on $E$.

\begin{definition} \label{canonicalconn} We say the map $\rho$ is the \emph{canonical connection} induced by the projective basis $B,B^*$.
We say $E$ is a $\Diff_\Z(A)$-module if the
  canonical map $\rho$ from Equation \ref{canonical} is a ring homomorphism.
\end{definition}

\begin{example} The free $A$-module of rank $n$.
\end{example}

Assume $E=A\{e_1,..,e_n\}$ is a free $A$-module of rank $n$ on the
basis $B=\{e_1,..,e_n\}$ and let $B^*\{x_1,..,x_n\}$ with $x_i=e_i^*\in E^*$. It follows 
$B,B^*$ is a projective basis for the finitely generated projective
$A$-module $E$. The map $\rho$ in this case is the canonical map
\[ \rho:\Diff^l_\Z(A)\rightarrow \Diff^l_\Z(E)\]
defined by
\[ \rho(\partial)(\sum_{i=1}^n a_ie_i)=\sum_{i=1}^n\partial(a_i)e_i.\]
Hence a differential operator $\partial\in \Diff^l_\Z(A)$ acts in each
coordinate for $E$. It follows the free $A$-module $E$ of rank $n$ is
a $\Diff_\Z(A)$-module in a canonical way.

The curvature of the $l$-connection $\nabla^l$ is related to the way
the map $\rho$ deviates from being a map of associative rings.

Recall the following definition from \cite{ogus} Proposition 2.10:
\begin{definition} A \emph{stratification} on $E$ is a collection of
  isomorphisms
\[ \eta_l:P^l\otimes_A E \rightarrow E\otimes_A P^l \]
such that
\begin{align}
&\label{str1}\text{$\eta_l$ is $P^l$-linear.}\\
&\label{str2}\text{$\eta_l$ and $\eta_k$ are compatible via
  restriction maps.}\\
&\label{str3}\text{$\eta_0$ is the identity map.}\\
&\label{str4}\text{The cocycle condition holds for $\eta_l$.}
\end{align}
\end{definition} 

See \cite{ogus} for a precise description of the \emph{cocycle
  condition}.
We say the system $\{\theta_l\}$ of right $A$-linear maps
\[ \theta_l:E\rightarrow E\otimes_A P^l \]
for all $l\geq 1$ with $\theta_0=id_E$ is an
\emph{$\infty$-connection} on $E$. 

Consider the following diagram:
\[
\label{curvature}
\diagram E\otimes_A P^{l+k} \rto^{id\otimes \delta^{l,k}} & E\otimes_A
P^l\otimes_A P^k \\
 E \uto^{\theta_{l+k}} \rto^{\theta_k} & E\otimes_A P^k
 \uto^{\theta_l\otimes id}
\enddiagram.
\]
Let $\phi_0^{l,k}=id\otimes \delta^{l,k}\circ \theta_{l+k}$ and
$\phi_1^{l,k}=\theta_l\otimes id \circ \theta_k$.

\begin{definition} \label{lkcurvature} Let $K^{l,k}=\phi_1^{l,k}-\phi_0^{l,k}$ be the
  $(l,k)$-curvature of the of the $\infty$-connection $\{\theta_l\}$.
We say the $\infty$-connection $\{\theta_l\}$ is \emph{flat} if
$K^{l,k}=0$ for all $l,k\geq 1$.
\end{definition}

\begin{proposition} The following data are equivalent:
\begin{align}
&\label{s1}\text{A stratification on $E$.}\\
&\label{s2}\text{A flat $\infty$-connection on $E$}.\\
&\label{s3}\text{A left $\Diff_\Z(A)$-module structure on $E$.}
\end{align}
\end{proposition}
\begin{proof} For a proof see \cite{ogus}, Proposition 2.11.
\end{proof}

\begin{example} Stratifications on finitely generated projective modules.\end{example}

In \cite{ogus}, 2.17 is is proved that if $A$ is a finitely generated
algebra over a field $K$ and $E$ is a finitely generated left $A$-module
with a stratification then $E$ is locally free, hence projective.
As indicated in Example \ref{TWOsphere}: The canonical connection
\[ \rho:\Diff_\Z(A)\rightarrow \Diff_\Z(E) \]
on a finitely generated projective module $E$
constructed using a projective basis $B,B^*$ is seldom a morphism of
rings. Hence the $A$-module $E$ seldom has a stratification in the sense
of \cite{ogus} induced by a projective basis. It is not clear if every
connection
\[ \nabla: E\rightarrow \Omega\otimes_A E \]
is induced by a projective basis $B,B^*$ for $E$. Given two
connections $\nabla, \nabla'$ their difference $\phi=\nabla-\nabla'$
is an $A$-linear map
\[ \phi: E\rightarrow \Omega\otimes_A E\]
and such a map may be constructed using a projective basis. It is thus
unlikely there is an action
\[ \rho:\Diff_K(A)\rightarrow \Diff_K(E) \]
where $\rho$ is a morphism of associative rings. Hence it is unlikely
a finitely generated projective $A$-module $E$ has a stratification.
The obstruction to the existence of a stratification is given by the
$(l,k)$-curvature $K^{(l,k)}$ from Definition \ref{lkcurvature}

\begin{example} Associative subrings of $\End_\Z(A)$.\end{example}

In general one may do the following: For any associative subring $R$
of $\End_\Z(A)$ containing $A$ we get a left action
\[ \eta:R\rightarrow \End_\Z(E)\]
defined by
\[\eta(\phi)(e)=\sum_{i=1}^n \phi(x_i(e))e_i.\]
If $a\in A\subseteq R$ acts via $\phi_a$ we get
\[ \eta(a)(e)=\sum_{i=1}^n \phi_a(x_i(e))e_i=\sum_{i=1}^n ax_i(e)e_i=\]
\[ \sum_{i=1}^nx_i(ae)e_i=ae \]
hence the $R$-structure on $E$ induced by $\eta$ extends the left
$A$-module structure on $E$. 

The map $\eta$ is not a map of associative rings in general.
If the map $\eta$ is a map of associative rings we  get a structure of left
$R$-module on $E$ extending the left $A$-module structure.
We say the $A$-module structure \emph{lifts} to $R$. This is a special case
of a general problem in
deformation theory: Given a map of rings $\phi:A\rightarrow B$ and a
left $A$-module $E$. One wants to extend this structure and define a
left $B$-module structure on $E$ restricting to the $A$-module structure via the map $\phi$.
As indicated in Example \ref{TWOsphere} for a finitely generated
projective $A$-module $E$ and a subring $R$ of $\End_\Z(A)$ the left
$A$-module structure on $E$ seldom lift to a left $R$-module structure. 

\textbf{Acknowledgements:} Thanks to Jean Fasel for references on the Bass-Quillen conjecture.

\end{document}